\documentclass[12pt]{amsart}
\usepackage{latexsym}
\usepackage{bbm}
\usepackage{amssymb,amsmath,mathrsfs, esint}
\usepackage{amsthm}
\usepackage{dsfont}

\newtheorem{theorem}{Theorem}

\theoremstyle{remark}
\newtheorem{rmk}{Remark}

\usepackage[hypertexnames=false, colorlinks, citecolor=black, linkcolor=blue, urlcolor=red]{hyperref}

\newcommand{\wt}{\widetilde}

\newcommand{\ci}[1]{_{ {}_{\scriptstyle #1}}}

\newcounter{vremennyj}

\numberwithin{equation}{section}

\begin{document}

\title[Generalized maximal operator]{A note on two weight bounds for the generalized Hardy-Littlewood Maximal operator}

\author{Amalia Culiuc}

\address{Department of Mathematics, Brown University, 151 Thayer
Str./Box 1917,      
 Providence, RI  02912, USA }
\email{amalia@math.brown.edu}



\begin{abstract}
We give a straighforward proof of the two weight estimates of the generalized maximal operator under Sawyer type testing conditions. The proof relies on the Martingale Carleson Embedding Theorem.
\end{abstract}

\maketitle

\bibliographystyle{plain}

\thispagestyle{empty}

\section{Introduction} 

Let $(X,\Sigma, \omega)$ be a $\sigma$-finite measure space with a filtration $\Sigma_n, n \in \mathbb{Z}$ (an increasing sequence of $\sigma$-algebras with $\Sigma_n\subset \Sigma$). We make the assumption that for each $\sigma$-algebra $\Sigma_n$ there exists a countable collection $\mathcal{D}_n$ of disjoint sets of finite measure such that any set in $\Sigma_n$ can be written as a union of sets in $\mathcal{D}_n$. The elements of any $\mathcal{D}_n$ are pairs of the form $(Q,n)$. However, in what follows, by slightly abusing notation, we will ignore the dependence on $n$ and simply write $Q$ instead.

Let $\mathcal{D}=\bigcup_n\mathcal{D}_n$. A classical example of such a filtration is that given by a dyadic lattice on $\mathbb{R}^n$ and in fact we may often refer to $\mathcal{D}$ as a lattice and to its elements as cubes. In spite of this language, we will not be making any further assumptions on the underlying structure of the space, including for example any assumptions about the homogeneity of $X$ with the measure $\omega$. 

Let 
$\mu$ and $\nu$ be measures, finite on all $Q\in\mathcal{D}$. 

For a sequence of functions $a=\{a\ci Q\}\ci{Q\in \mathcal{D}}$, $a\ci Q: X\rightarrow [0,\infty)$ indexed by the sequence of dyadic cubes, define  the operator $M_a^q$, given by
\begin{equation*}
M_a^q f\mu (x)=\left(\sum_{\substack{Q\in \mathcal{D}\\ x\in Q}}\left|\left(\int_Q f d\mu \right) a\ci Q(x) \mathbbm{1}\ci Q(x)\right|^q\right)^{1/q} \text{  for $1< q<\infty$}
\end{equation*}
and  
\begin{equation*}
M^\infty_a f\mu(x) = \sup \left\{ \left|\left(\int_Q f d\mu \right) a\ci Q(x) \mathbbm{1}\ci Q(x)\right| \,:\, Q\in\mathcal{D} \quad\text{such that } x\in Q \right\}\,.
\end{equation*}


In what follows we will show that under a so-called Sawyer type testing condition, the operator $f\mapsto M_a^q f\mu$ is bounded $L^p(\mu)\rightarrow L^p(\nu)$ for $p\leq q$. Such conditions were named after E.~Sawyer, who introduced them in \cite{Sawyer} when studying the two weight estimates for  the classical maximal operator $M$. The testing condition presented in \cite{Sawyer} essentially amounts to testing the uniform estimates 
on characteristic functions of dyadic cubes. Later, in \cite{Sawyer2}, Sawyer proved that for operators such as fractional integrals, Poisson kernels, and other nonnegative kernels, the two weight estimate still holds if one assumes the testing condition not only on the operator itself, but also on its formal adjoint. For the positive martingale operators such results were obtained in \cite{ntv99} ($p=2$) and later in \cite{lsut} (general $p$), see  also \cite{Treil} for a simpler proof.

The main result of this paper is a simple proof of the theorem below. Formally, this result is new, because we allow $a\ci Q$ to be nonnegative functions and we also consider the case $q<\infty$.

For $Q\in\mathcal{D}$ define the truncated operator $M\ci{a,Q}^q$, 
\[
M\ci{a,Q}^q f\mu (x):= \left(\sum_{\substack{R\in \mathcal{D},R\subset Q\\ x\in R}}\left|\left(\int_R f d\mu \right) a\ci R(x) \mathbbm{1}\ci R(x)\right|^q\right)^{1/q} , 
\]
with the obvious modification for $q=\infty$. 
\begin{theorem}\label{main}
Let $1<p\leq q\leq \infty$. The operator $M_{a}^q$ 
 satisfies
\begin{align}
\label{MainEst}
\| M_a^q f\mu\|\ci{L^p(\nu)}^p \le A \|f\|\ci{L^p(\mu)} \qquad \forall f\in L^p(\mu)
\end{align}
if and only if the truncation $M\ci{a,Q}^q$ satisfies the following testing condition 
\begin{equation}\label{testing}
\|M_{a,Q}^q(\mathbbm{1}\ci Q \mu)\|_{L^p(\nu)}\leq B \mu(Q)^{1/p}, \text{  for any } Q\in \mathcal{D}. 
\end{equation}
Moreover, for the best constants $A$ and $B$ we have $B\le A\leq C(p) B$, 
\[
C(p) = \left(\left(1+{1}/{p}\right)^{p+1}p\right)^{1/p}p', 
\]
where $p'$ is the H\"{o}lder conjugate of $p$,  $1/p+ 1/{p'}=1$.  
\end{theorem}

For $p=q=\infty$, the result is trivial with $A=B$. Note that $\lim_{p\to\infty} C(p) =1$. 
\begin{rmk}
The  classical dyadic (martingale) maximal operator $M$ is a particular case of operator $M^q_a$ where $q=\infty$,  $a\ci Q \equiv \omega(I)^{-1} \mathbbm{1}\ci Q$,  and $\mathcal D$ is a dyadic lattice in $\mathbb R^n$. Therefore, one can view the $M_a^q$ as a generalization of the classical martingale maximal function.

In \cite{Sawyer}, E.~Sawyer considered slightly more general maximal operators $M=M^\alpha$ which are a particular case of our $M_a^q$ with $q=\infty$ and   $a\ci Q \equiv \omega(I)^{-\alpha} \mathbbm{1}\ci Q$, $0<\alpha\le1$.  He characterized the measures $\wt\mu$, $\nu$ and $\omega$ for which the inequality
\begin{align}
\label{SawyerEst}
\| Mf\omega\|\ci{L^p(\nu)} \le A \|f\|\ci{L^p(\wt\mu)}\qquad \forall f\in L^p(\wt\mu)
\end{align}
holds. 

Note that without loss of generality one can assume that $\wt\mu$ is absolutely continuous with respect to $\omega$, $d\wt\mu = wd\omega$ (adding a singular part to $\wt\mu$ does not change \eqref{SawyerEst}). 
So,  making the standard change of weight  $f\mapsto w^{p'/p} f$ and denoting $\mu := w^{-p'/p} \omega$ we transform  the above estimate \eqref{SawyerEst} to \eqref{MainEst}. Then in this notation the necessary and sufficient condition obtained by E.~Sawyer is exactly the testing condition \eqref{testing}.

%
%
%
%

For the classical dyadic maximal operator $M=M^\alpha$, the truncation $M_{a, Q}^\infty \mathbbm{1}\ci{Q} \mu$ defined above is equivalent to $\mathbbm{1}\ci Q M\mathbbm{1}\ci{Q}\mu$, so up to a change of measure, our setup is identical to \cite{Sawyer}. 
\end{rmk}

\begin{rmk}
The reduction to the two measure setup, eliminating the underlying measure $\omega$, is now considered standard for weighted estimates. In the three measures setup for the classical maximal operator as in \cite{Sawyer} 
all the information about $\omega$ is captured by the (constant in this case) functions  $a\ci Q$. 

To obtain Sawyer's estimate for the non-martingale maximal function, one can use the two weighted estimate for the dyadic case and proceed by an averaging argument. This reasoning is pretty standard and we will not discuss it at all  in this paper.

Our proof simplifies the argument in \cite{Sawyer} and gives a formally stronger result: in particular, the coefficients $a\ci Q$ do not need to be constants.  It also has the additional benefit of placing the Hardy-Littlewood maximal function in the context of a wide range of similar operators.
\end{rmk}

The proof we present relies on the stopping time construction presented in \cite{Treil} and the Martingale Carleson Embedding Theorem stated below.

Denote
\[
\fint_Q f d\mu=\frac{1}{\mu(Q)}\int_Q f d\mu. 
\]

\begin{theorem}{(Martingale Carleson Embedding Theorem)}\label{Carleson}
Let $\mu$ be a measure on $X$ and let $\{w\ci Q\}\ci{Q\in \mathcal{D}}$, $w\ci Q\geq 0$ be a sequence satisfying the following condition: 
\begin{equation*}
\sum_{Q\subset R, R\in \mathcal{D}}w\ci Q\leq A\mu(R), \, \text{ for any cube $R\in \mathcal{D}$}.
\end{equation*}
Then for any measurable function $f\geq 0$ and for any $p\in(1,\infty)$,
\begin{equation*}
\sum_{Q\in \mathcal{D}}\left(\fint_Q f d\mu\right)^p w\ci Q\leq (p')^p  A\|f\|^p_{L^p(\mu)}.
\end{equation*}
\end{theorem}

The Carleson Embedding Theorem with the constant $(p')^p$ can be proved as a straightforward consequence of the one weight $L^p$ boundedness of the classical Hardy-Littlewood maximal function (see \cite{Treil}). Other arguments that include the sharp constant have been given by Nazarov, Treil, and Volberg \cite{ntv} for $p=2$ and Lai \cite{Lai} for $p\ne2$ using Bellman function techinques. The exact Bellman function for $p>1$ was originally computed by Melas in \cite{melas}, but the sharp constant was not explicitly stated. 

\subsection*{Acknowledgement}
The author would like to thank her advisor, Prof. Sergei Treil, for suggesting the idea and for all the guidance in the process of writing this paper.

\section{Proof of the main result}
We aim to prove Theorem \ref{main}. Again, we make the remark that while we may refer to the elements $Q\in \mathcal{D}$ as cubes, $\mathcal{D}$ is not necessarily assumed to be the dyadic lattice on $\mathbb{R}^n$.

\begin{proof}
The necessity of the testing condition and the estimate $B\leq A$ are trivial: if $M_a^q$ is bounded on $L^p(\mu)$ functions, it is, in particular, bounded on caracteristic functions. Thus, by testing $M_a^q$ on the functions $\mathbbm{1}_Q$, we obtain condition \eqref{testing}. 

To prove sufficiency, we begin by constructing a collection of stopping cubes $\mathcal{G}\subset \mathcal{D}$, following the definitions and notation in \cite{Treil}.

For any cube $Q\in \mathcal{D}$, define $\mathcal{D}(Q)$ to be the collection of subcubes of $Q$ in $\mathcal{D}$. For a fixed $r>1$, let $\mathcal{G}^*(Q)$ be the set of stopping cubes of $Q$, that is,
\begin{equation*}
\mathcal{G}^{*}(Q)=\left\{ R\in \mathcal{D}(Q),  R \text{ maximal, } \frac{1}{\mu(R)}\int_R f d\mu\geq r\frac{1}{\mu(Q)}\int_Q f d\mu \right\},
\end{equation*}
where maximality is considered with respect to the partial ordering given by inclusion. 

Denote by $\mathcal{E}(Q)$ the collection of descendants of $Q$ that are not stopping cubes or descendants of the stopping cubes:
\begin{equation*}
\mathcal{E}(Q)=\mathcal{D}(Q)\setminus \bigcup_{P\in\mathcal{G}^{*}(Q)}\mathcal{D}(P).
\end{equation*}
Note that, by definition, for any $R\in \mathcal{E}(Q)$,
\begin{equation}\label{bound}
\fint_R f d\mu< r \fint_Q f d\mu.
\end{equation}
Also note that
\begin{equation}\label{condition}
\sum_{R\in\mathcal{G}^{*}(Q)} \mu(R) = \mu\left(\bigcup_{R\in\mathcal{G}^{*}(Q)}R\right)\leq \frac{\mu(Q)}{r}.
\end{equation}

To construct the collection $\mathcal{G}$ of stopping cubes, let $N$ be a fixed large positive integer and define the first generation $\mathcal G_1$ as
\begin{equation*}
\mathcal{G}_1= \mathcal{D}_{-N}.
\end{equation*}
Then, to obtain the subsequent generations, apply the inductive formula 
\begin{equation*}
\mathcal{G}_{n+1}=\bigcup_{Q\in \mathcal{G}_n}\mathcal{G}^{*}(Q).
\end{equation*}
Define the collection of stopping cubes $\mathcal{G}$ to be the union 
\begin{equation*}
\mathcal{G}=\bigcup_{n=1}^\infty \mathcal{G}_n.
\end{equation*}
Equation \eqref{condition} implies that 
\begin{equation}\label{measurecond}
\sum_{R\in \mathcal{G},R\subset Q}\mu(R)\leq \frac{r}{r-1}\mu(Q), \, \forall \, Q\in \mathcal{D}.
\end{equation}

To prove the theorem it is sufficient to prove the uniform in $N$ bounds for the operator $M\ci{a}^{q,N}$, 
\[
M\ci{a}^{q,N} f:=\left( \sum_{n\ge-N}\sum_{Q\in\mathcal D_n} \left(\left(\fint_R f d\mu \right) \mu(R) a\ci R(x)\mathbbm{1}\ci R (x)\right)^q\right)^{1/q}  
\]
(with the obvious change for $q=\infty$) and then let $N\to\infty$. 

Given the construction of stopping moments, it is easy to see that
\begin{equation*}
\bigcup_{n=-N}^\infty\mathcal{D}\ci{n}=\bigcup_{Q\in \mathcal{G}}\mathcal{E}( Q ).
\end{equation*}
and that the sets $\mathcal E(Q)$ are disjoint. 

In the proof below we use notation for $1<q<\infty$. The proof for $q=\infty$ is absolutely the same (up to obvious changes in the notation).

Denoting
\begin{equation*}
F\ci Q(x)=\left(\sum_{R\in \mathcal{E}(Q)}\left(\left(\fint_R f d\mu \right)\mu(R) a\ci R(x)\mathbbm{1}\ci R (x)\right)^q\right)^{1/q},
\end{equation*}
we can write
$\displaystyle
M\ci{a}^{q,N} f = \left(\sum_{Q\in \mathcal{G}} F\ci Q^q\right)^{1/q}
$, 
so the proof amounts to bounding $\left \| \left(\displaystyle\sum_{Q\in\mathcal{G}} F\ci Q^q\right)^{1/q} \right\|_{L^p(\nu)}$. 

Since $\|x\|_{\ell^q} \le \|x\|_{\ell^p}$ for $q\ge p$, we can estimate
\begin{align}\label{est1}
\left \| \left(\sum_{Q\in\mathcal{G}} F\ci Q^q\right)^{1/q} \right\|_{L^p(\nu)}&\leq \left \| \left(\sum_{Q\in\mathcal{G}} F\ci Q^p\right)^{1/p} \right\|_{L^p(\nu)}\\
&=\left( \int \sum_{Q\in\mathcal{G}} F\ci Q^p d\nu\right)^{1/p}\notag\\
&=\left(\sum_{Q\in\mathcal{G}} \int F\ci Q^p d\nu\right)^{1/p}\notag\\
\notag
&=\left(\sum_{Q\in\mathcal{G}} \| F\ci Q\|^p_{L^p(\nu)}\right)^{1/p}.
\end{align}
By definition,
\begin{align}\label{est2}
\left\| F\ci Q\right\|^p_{L^p(\nu)}&=\left\|\left(\sum_{R\in \mathcal{E}(Q)}\left(\left(\fint_R f d\mu \right)\mu(R) a\ci R\mathbbm{1}\ci R\right)^q\right)^{1/q} \right\|^p_{L^p(\nu)}
\\
&\leq \left\|\left(\sum_{R\in \mathcal{E}(Q)}\left(\left(r \fint_Q f d\mu \right)\mu(R) a\ci R\mathbbm{1}\ci R \right)^q\right)^{1/q} \right\|^p_{L^p(\nu)} &&\text{ by  \eqref{bound}}\notag\\
&=r^p \left(\fint_Q f d\mu\right)^p \left\|\left(\sum_{R\in \mathcal{E}(Q)}\left(  
\left( \int_R \mathbbm{1}\ci{Q} d\mu \right)a\ci R \mathbbm{1}\ci R\right)^q
\right)^{1/q}\right\|^p_{L^p(\nu)} \notag \\
\notag
& \le r^p B^p \left(\fint_Q  f d\mu\right)^p \mu(Q) && \text{by \eqref{testing} } .
\end{align}

Therefore, from inequalities \eqref{est1} and \eqref{est2} we obtain
\begin{align*}
\left \| \left(\sum_{Q\in\mathcal{G}} F\ci Q^q\right)^{1/q} \right\|_{L^p(\nu)}& \leq \left(\sum_{Q\in\mathcal{G}} \| F\ci Q\|^p_{L^p(\nu)}\right)^{1/p}\\
&\leq r B \left(\sum_{Q\in\mathcal{G}} \left(\fint_Q d\mu\right)^p \mu(Q)\right)^{1/p}.
\end{align*}

Apply Theorem \ref{Carleson}, taking
\begin{equation*}
w\ci Q = \left\{
     \begin{array}{lr}
       \mu(Q) & : Q \in \mathcal{G}\\
       0 & : Q \notin \mathcal{G}
     \end{array}
   \right.
\end{equation*}

Equation \eqref{measurecond} shows that the sequence $w\ci Q$ satisfies the Carleson measure condition. Hence
\[
\sum_{Q\in\mathcal{G}} \left(\fint_Q f d\mu\right)^p \mu(Q)\leq \frac{r}{r-1}(p')^p \|f\|^p_{L^p(\mu)}.
\]

Consequently,

\begin{equation*}
\left\|M\ci{a}^{q,N}f\mu\right\|^p_{L^p(\nu)}\leq \frac{r^{p+1}}{r-1}(p')^pB^p\|f\|_{L^p(\mu)}^p.
\end{equation*}
In particular, since no assumption was made on $r$ other than $r>1$, one can consider the minimal value of the right hand side constant, which is attained at $\frac{p+1}{p}$. Then
\begin{equation*}
\left\|M\ci {a}^{q,N}f\mu\right\|^p_{L^p(\nu)}\leq \left(1+\frac{1}{p}\right)^{p+1}p(p')^pB^p\|f\|_{L^p(\mu)}^p.
\end{equation*}

The right hand side above does not depend on the choice of $N$. Taking the limit as $N$ approaches $\infty$ completes the proof for $M_a^q$.
\end{proof}


\begin{thebibliography}{5}
\bibitem{lsut} M. Lacey, E. Sawyer, I. Uriarte-Tuero, \textit{Two weight inequalities for discrete positive operators}, (2009) arXiv:0911.3437. 
\bibitem{Lai} J. Lai, \textit{Two weight problems and Bellman functions on filtered  spaces}, (2015), PhD thesis, Brown University.
\bibitem{melas} A. Melas, \textit{The Bellman functions of dyadic-like maximal operators and related inequalities}, Adv. Mat., \textbf{192},(2005), no. 2, 310--340.

\bibitem{ntv99} F. Nazarov, S. Treil, A. Volberg, \emph{The Bellman function and two weight inequalities for Haar 
multipliers}, (with F. Nazarov and A.  Volberg), Journal of the 
American Mathematical Society, {\bf 12}(1999), No 4, 909--928.

\bibitem{ntv} F. Nazarov, S. Treil, A. Volberg, \textit{Bellman functions in stochastic control and harmonic analysis}, Systems, Approximation, Singular Integral Operators, and Related Topics (Bordeaux, 2000), 393--423, Oper.Theory Adv. Appl., \textbf{129}, Birkh\"{a}user Basel, 2001
\bibitem{Sawyer} E. Sawyer, \textit{A characterization of a two-weight norm inequality for maximal operators}, Studia Math., \textbf{75} (1982), no.1, 1--11.
\bibitem{Sawyer2} E. Sawyer, \textit{A characterization of a two-weight norm inequality for fractional and Poisson integrals}, Trans. Amer. Math. Soc., \textbf{308} (1998), no.2, 533--545.
\bibitem{Treil} S. Treil, \textit{A remark on two weight estimates for positive dyadic operators}, in:  K.~Gröchenig, Yu.~Lyubarskii, K.~Seip (Editors), \emph{Operator-Related Function Theory and Time-Frequency Analysis: The Abel Symposium 2012 (Abel Symposia 9),} Springer, 2014, p.~185--195; see also arXiv:1201.1455. 
\end{thebibliography}
\end{document}